\documentclass[a4paper,twoside]{article}
\usepackage{amssymb}

\usepackage{makeidx}
\usepackage{amsbsy}
\usepackage{amsmath}
\usepackage{theorem}
\usepackage{graphicx}
\usepackage{epsfig}

\evensidemargin\oddsidemargin
\numberwithin{equation}{section}
\newtheorem{theorem}{Theorem}[section]
\newtheorem{corollary}[theorem]{Corollary}

\newtheorem{proposition}[theorem]{Proposition}

{\theorembodyfont{\rmfamily} 

}
{\theorembodyfont{\rmfamily} 
\newtheorem{remark}[theorem]{Remark}
}
{\theorembodyfont{\rmfamily} 

}
{\theorembodyfont{\rmfamily} 
\newtheorem{example}[theorem]{Example}
}
{\theorembodyfont{\rmfamily} 

}
{\theorembodyfont{\rmfamily} 
\newtheorem{final remarks and questions}[theorem]{Final remarks and questions}
}
\begin{document}

\title{\textbf{Crepant Resolutions of Gorenstein Toric Singularities and Upper
Bound Theorem}}
\author{\thanks{{\scriptsize GSRT-fellow supported by the European Union and the
Greek Ministry of Research and Technology.}\medskip \newline
} \textbf{Dimitrios I. Dais} \\
\noindent {\scriptsize Mathematics Department, Section of Algebra and
Geometry, University of Ioannina}\\
{\scriptsize GR-45110 Ioannina, Greece. E-mail: ddais@cc.uoi.gr}}
\date{}
\maketitle

\begin{abstract}
\noindent A necessary condition for the existence of torus-equivariant
crepant resolutions of Gorenstein toric singularities can be derived by
making use of a variant of the classical Upper Bound Theorem which is valid
for simplicial balls.
\end{abstract}

\section{Introduction\label{SECTION1}}

\noindent {}Let $d$ be a positive integer, $\sigma \subset \mathbb{R}^{d+1}$
a rational, $\left( d+1\right) $-dimensional strongly convex polyhedral cone
(w.r.t. the lattice $\mathbb{Z}^{d+1}$), and 
\begin{equation*}
U_{\sigma }=\text{Spec}\left( \mathbb{C}\left[ \sigma ^{\vee }\cap (\mathbb{Z%
}^{d+1})^{\vee }\right] \right)
\end{equation*}
the associated affine toric variety, where $\sigma ^{\vee }$ denotes the
dual of $\sigma $. (For the usual notions of toric geometry, see \cite{Oda}). 
As it is known (see e.g. \cite[\S 6]{Stanley1A}%
):

\begin{theorem}
\label{GORC}$U_{\sigma }$ is Gorenstein if and only if the set \emph{Gen}$%
\left( \sigma \right) $ of the minimal generating integral vectors of $%
\sigma $ lies on a ``primitive'' affine hyperplane, i.e., iff 
\begin{equation*}
\text{\emph{Gen}}\left( \sigma \right) \subset \mathbf{H}_{\sigma }=\left\{ 
\mathbf{x}\in \mathbb{R}^{d+1}\ \left| \ \left\langle \mathbf{m}_{\sigma },%
\mathbf{x}\right\rangle =1\right. \right\} ,
\end{equation*}
where $\mathbf{m}_{\sigma }\in (\mathbb{Z}^{d+1})^{\vee }$ is a primitive
vector belonging to the dual lattice.
\end{theorem}

\begin{remark}
(i) In this case, $\sigma $ supports the $d$-dimensional lattice polytope 
\begin{equation}
P_{\sigma }=\left\{ \mathbf{x}\in \sigma \ \left| \ \left\langle \mathbf{m}%
_{\sigma },\mathbf{x}\right\rangle =1\right. \right\} \subset \mathbf{H}%
_{\sigma }\cong \mathbb{R}^{d}  \label{PISIGMA}
\end{equation}
(w.r.t. the lattice $\mathbf{H}_{\sigma }\cap \mathbb{Z}^{d+1}\cong \mathbb{Z%
}^{d}$).\smallskip\ \newline
(ii) In fact, \textit{every} lattice $d$-polytope $P\subset \mathbb{R}^{d}$
can be considered as supported by a cone 
\begin{equation*}
\sigma _{P}=\left\{ \left( r,r\mathbf{x}\right) \in \mathbb{R}\oplus \mathbb{%
R}^{d}\ \left| \ \mathbf{x}\in P,\ r\in \mathbb{R}_{\geq 0}\right. \right\}
\subset \mathbb{R}^{d+1}
\end{equation*}
so that $U_{\sigma _{P}}$ is Gorenstein.
\end{remark}

\noindent {}The $\left( d+1\right) $-dimensional Gorenstein toric
singularities\footnote{%
Without loss of generality, we may henceforth assume that the cones $\sigma
\subset \mathbb{R}^{d+1}$ are $\left( d+1\right) $-dimensional, and that the
singularities under consideration have maximal splitting codimension. (The
orbit orb$\left( \sigma \right) \in U_{\sigma }$ is the unique fixed closed
point under the usual torus-action on $U_{\sigma }$.)} $\left( U_{\sigma },%
\text{orb}\left( \sigma \right) \right) $ constructed by cones $\sigma $
which support lattice $d$-polytopes $P=P_{\sigma }$ are to be subdivided
into three distinct classes\footnote{%
A lattice polytope $P$\ is called \textit{elementary} if the lattice points
belonging to it are exactly its vertices. A lattice simplex is said to be 
\textit{basic} (or \textit{unimodular}) if its vertices constitute a part of
a\emph{\ }$\mathbb{Z}$-basis of the reference lattice (or equivalently, if
its relative, normalized volume equals $1$). A lattice triangulation\emph{\ }%
$\mathcal{T}$ $\ $of a lattice polytope $P$ is defined to be \textit{maximal}
(resp. \textit{basic}), if it consists only of elementary (resp. basic)
simplices.}:\smallskip

\noindent {}\textbf{(A)} \textit{Terminal singularities} (whenever $P$ is an
elementary polytope but not a basic simplex).\smallskip\ 

\noindent {}\textbf{(B)} \textit{Canonical, non-terminal singularities which
do not admit any crepant resolution} (i.e., for which $P$ is a non-elementary 
polytope having no basic triangulations).\smallskip

\noindent {}\textbf{(C)} \textit{Canonical, non-terminal singularities
admitting crepant resolutions} (i.e., for which $P$ is a non-elementary
polytope possessing \textit{at least one}
basic triangulation).\smallskip

\noindent {}\textbf{Comments.} A complete classification of the members of
class \textbf{(A)} (up to analytic isomorphism) is obviously equivalent to
the classification of elementary polytopes (up to lattice automorphism). For
constructions of several families belonging to \textbf{(C)}, the reader is
referred to \cite{DHaZ}, \cite{DHH}, \cite{DH}, \cite{DHZ1}, \cite{DHZ2}. In
fact, for one- or two-parameter Gorenstein cyclic quotient singularities, it
is possible to decide definitely if they belong to class \textbf{(A)}, 
\textbf{(B)} or \textbf{(C)}, by just checking some concrete
number-theoretic (necessary and sufficient existence-) conditions (see \cite
{DH} and \cite{DHH}, respectively). On the other hand, \ for \textit{general}
Gorenstein toric (not necessarily quotient-) singularities, a \textit{%
necessary} \textit{condition} for the existence of crepant resolutions can
be created via an UBT for simplicial balls, as we shall see below in Thm. 
\ref{MAIN}. Hence, its ``violation'' may be used to produce families of such
singularities belonging to \textbf{(B)}.

\section{Basic facts about UBT's\label{SECTION2}}

\noindent {}\textsc{Notation}. (i) The $\mathbf{f}$-\textit{vector }$\mathbf{%
f}\left( \mathcal{S}\right) =\left( \frak{f}_{0}\left( \mathcal{S}\right) ,%
\frak{f}_{1}\left( \mathcal{S}\right) ,\ldots ,\frak{f}_{d-1}\left( \mathcal{%
S}\right) \right) $ of a polyhedral $\left( d-1\right) $-complex\emph{\ }$%
\mathcal{S}$ is defined by setting for all $i$, $0\leq i\leq d-1$, 
\begin{equation*}
\frak{f}_{i}\left( \mathcal{S}\right) :=\#\ \left\{ i\text{-dimensional
faces of\emph{\ }}\mathcal{S}\right\}
\end{equation*}
(under the usual conventional extension: $\frak{f}_{-1}\left( \mathcal{S}%
\right) :=1$). The\emph{\ }coordinates of the $\mathbf{h}$-\textit{vector } $%
\mathbf{h}\left( \mathcal{S}\right) =\left( \frak{h}_{0}\left( \mathcal{S}%
\right) ,\frak{h}_{1}\left( \mathcal{S}\right) ,\ldots ,\frak{h}_{d-1}\left( 
\mathcal{S}\right) ,\frak{h}_{d}\left( \mathcal{S}\right) \right) $ of such
an $\mathcal{S}$ are defined by the equations 
\begin{equation}
\frak{h}_{j}\left( \mathcal{S}\right) =\sum_{i=0}^{j}\ \left( -1\right)
^{j-i}\ \tbinom{d-i}{d-j}\ \frak{f}_{i-1}\left( \mathcal{S}\right) .
\label{HD-REL}
\end{equation}

\noindent {}(ii) For a $d$-dimensional polytope $P,$ the \textit{boundary
complex }$\mathcal{S}_{\partial P}$ of $P$ is defined to be the $\left(
d-1\right) $-dimensional polyhedral complex consisting of the proper faces
of $P$ together with $\varnothing $ and having $\partial P$ as its support. $%
\mathcal{S}_{\partial P}$ is a polyhedral $\left( d-1\right) $-sphere. $%
\mathcal{S}_{\partial P}$ is a geometric pure simplicial complex (in fact, a
simplicial $\left( d-1\right) $-sphere) if and only if $P$ is a simplicial
polytope. The $\mathbf{f}$-vector $\mathbf{f}\left( P\right) $ of a $d$%
-polytope $P$\emph{\ }is by definition\emph{\ }the $\mathbf{f}$-vector $%
\mathbf{f}\left( \mathcal{S}_{\partial P}\right) $ of its boundary
complex.\medskip

\noindent {}(iii) We denote by CycP$_{d}\left( k\right) $ the \textit{cyclic}
$d$-\textit{polytope} with $k$ vertices. As it is known, the number of its
facets equals 
\begin{equation}
\frak{f}_{d-1}\left( \text{CycP}_{d}\left( k\right) \right) =\tbinom{%
k-\left\lceil \frac{d}{2}\right\rceil }{\left\lfloor \frac{d}{2}%
\right\rfloor }+\tbinom{k-1-\left\lceil \frac{d-1}{2}\right\rceil }{%
\left\lfloor \frac{d-1}{2}\right\rfloor }  \label{FDCYCLIC}
\end{equation}

\noindent This follows from Gale's evenness condition and the fact that CycP$%
_{d}\left( k\right) $ is $\left\lfloor \frac{d}{2}\right\rfloor $%
-neighbourly (cf. \cite[p. 24]{Ziegler}).\medskip\ 

\noindent {}(iv) Classical UB and LB-Theorems for simplicial spheres (see 
\cite{Stanley1} and \cite{KALAI}):

\begin{theorem}[Stanley's Upper Bound Theorem for Simplicial Spheres]
\label{UBTSS}\negthinspace \negthinspace\ \newline
The $\mathbf{f}$-vector coordinates of a simplicial $\left( d-1\right) $%
-sphere $\mathcal{S}$ with $\frak{f}_{0}\left( \mathcal{S}\right) =k$
vertices satisfy the following inequalities\emph{:} 
\begin{equation*}
\frak{f}_{i}\left( \mathcal{S}\right) \leq \frak{f}_{i}\left( \emph{CycP}%
_{d}\left( k\right) \right) \ ,\ \forall \ i\ ,\ 0\leq i\leq d-1.
\end{equation*}
\end{theorem}

\begin{theorem}[Lower Bound Theorem for Simplicial Spheres]
\label{LBTH}The $\mathbf{h}$-vector coordinates of a simplicial $\left(
d-1\right) $-sphere $\mathcal{S}$ with $\frak{f}_{0}\left( \mathcal{S}%
\right) =k$ vertices satisfy the following inequalities\emph{:} 
\begin{equation*}
\frak{h}_{1}\left( P\right) =k-d\leq \frak{h}_{i}\left( P\right) \ ,\
\forall \ i\ ,\ 2\leq i\leq d.
\end{equation*}
\end{theorem}

\noindent {}Besides them we need certain variants for simplicial \textit{%
balls}.

\begin{proposition}[``$\mathbf{h}$ of $\partial "-$Lemma]
\label{STANLEM}Let $\mathcal{S}$ be a $d$-dimensional Cohen-Macaulay closed
pseudomanifold with non-empty boundary $\partial \mathcal{S}$. Then 
\begin{equation}
\frak{h}_{i-1}\left( \partial \mathcal{S}\right) -\frak{h}_{i}\left(
\partial \mathcal{S}\right) =\frak{h}_{\left( d+1\right) -i}\left( \mathcal{S%
}\right) -\frak{h}_{i}\left( \mathcal{S}\right) ,\ \ \forall i,\ 0\leq i\leq
d+1\smallskip  \label{STAL2}
\end{equation}
\emph{(}under the convention\emph{: }$\frak{h}_{-1}\left( \partial \mathcal{S%
}\right) =0$\emph{).}
\end{proposition}

\noindent \textit{Proof. }See Stanley (\cite[2.3]{Stanley5}).\hfill$\square
\bigskip $

\noindent {}Working with \textit{Buchsbaum complexes}, Schenzel \cite
{Schenzel} proved the following:

\begin{theorem}[Schenzel's Upper Bound Theorem]
Let $\mathcal{S}$ be a $d$-dimensional Buchsbaum complex\emph{\footnote{%
A simplicial complex $\mathcal{S}$ is a Buchsbaum complex over a field $%
\mathbf{k}$ if and only if it is pure and the localizations $\mathbf{k}\left[
\mathcal{S}\right] _{\wp }$ of $\mathbf{k}\left[ \mathcal{S}\right] $ w.r.t.
prime ideals $\wp \neq \mathbf{k}\left[ \mathcal{S}\right]
_{+}(=\bigoplus_{\nu >0}\left( \mathbf{k}\left[ \mathcal{S}\right] _{\nu
}\right) $ are Cohen-Macaulay. (For instance, homology $d$-manifolds without
boundary, or homology $d$-manifolds whose boundary is a homology $\left(
d-1\right) $-manifold without boundary, are Buchsbaum). Moreover, $\mathcal{S%
}$ is Cohen-Macaulay over $\mathbf{k}$ if an only if $\mathcal{S}$ is
Buchsbaum over $\mathbf{k}$ and dim$_{\mathbf{k}}\widetilde{H}%
_{j}^{\,}\left( \mathcal{S};\mathbf{k}\right) =0$, for all $i,0\leq i\leq
d-1 $, while dim$_{\mathbf{k}}\widetilde{H}_{d}^{\,}\left( \mathcal{S};%
\mathbf{k}\right) =\left( -1\right) ^{d}\ \widetilde{\chi }\left( \mathcal{S}%
\right) $, with $\widetilde{\chi }\left( \mathcal{S}\right) $ the reduced
Euler characteristic.}} having $\frak{f}_{0}\left( \mathcal{S}\right) =\frak{%
b}$ vertices. Then for all $i$, $0\leq i\leq d+1$, the $\mathbf{h}$-vector
coordinates of $\mathcal{S}$ satisfy the inequalities \emph{\ } 
\begin{equation}
\frak{h}_{i}\left( \mathcal{S}\right) \leq \tbinom{\frak{b}-d+i-2}{i}-\left(
-1\right) ^{i}\ \tbinom{d+1}{i}\ \sum_{j=-1}^{i-2}\ \left( -1\right) ^{j}\ 
\text{\emph{dim}}_{\mathbf{k}}\widetilde{H}_{j}^{\,}\left( \mathcal{S};%
\mathbf{k}\right)  \label{HSCHEN}
\end{equation}
\emph{(}where $\widetilde{H}_{j}^{\,}\left( \mathcal{S};\mathbf{k}\right) $
are the reduced homology groups of $\mathcal{S}$ with coefficients in a
field $\mathbf{k}.\emph{)}$
\end{theorem}

\begin{corollary}
Let $\mathcal{S}$ denote a simplicial $d$-dimensional $\mathbf{ball}$ with $%
\frak{f}_{0}\left( \mathcal{S}\right) =\frak{b}$ vertices. Then for all $%
i,\,\,\,0\leq i\leq d$, the $\mathbf{f}$-vector of $\mathcal{S}$ satisfies
the following inequalities\emph{:\medskip } 
\begin{equation}
\frak{f}_{i}\left( \mathcal{S}\right) \leq \frak{f}_{i}\left( \emph{CycP}%
_{d+1}\left( \frak{b}\right) \right) -\sum\limits_{j=d-i}^{\left\lfloor 
\frac{d}{2}\right\rfloor \smallskip }\ \tbinom{j}{d-i}\ \left( \frak{h}%
_{j}\left( \partial \mathcal{S}\right) -\frak{h}_{j-1}\left( \partial 
\mathcal{S}\right) \right) ,  \label{WEAKUBTSB}
\end{equation}
\end{corollary}

\noindent \textit{Proof. }Introduce the auxiliary vector $\widetilde{\mathbf{%
h}}\left( \mathcal{S}\right) =\left( \widetilde{\frak{h}}_{0}\left( \mathcal{%
S}\right) ,\ldots ,\widetilde{\frak{h}}_{d+1}\left( \mathcal{S}\right)
\right) $ with{\small 
\begin{equation*}
\widetilde{\frak{h}}_{i}\left( \mathcal{S}\right) :=\left\{ 
\begin{array}{lll}
\frak{h}_{i}\left( \mathcal{S}\right) & , & \text{for \ \ \ }0\leq i\leq
\left\lfloor \frac{d+1}{2}\right\rfloor \\ 
\  &  &  \\ 
\frak{h}_{i}\left( \mathcal{S}\right) -\left( \frak{h}_{d-i}\left( \partial 
\mathcal{S}\right) -\frak{h}_{d+1-i}\left( \partial \mathcal{S}\right)
\right) & , & \text{for \ \ \ }\left\lfloor \frac{d+1}{2}\right\rfloor
+1\leq i\leq d+1
\end{array}
\right.
\end{equation*}
}Since $\mathcal{S}$ is Cohen-Macaulay, $\mathcal{S}$ is a Buchsbaum
complex. Moreover, all reduced homology groups $\widetilde{H}_{j}^{\,}\left( 
\mathcal{S};\mathbf{k}\right) $ are trivial, which means that 
\begin{equation*}
\frak{h}_{i}\left( \mathcal{S}\right) \leq \frak{h}_{i}\left( \text{CycP}%
_{d+1}\left( \frak{b}\right) \right) =\tbinom{\frak{b}-d+i-2}{i}%
,\,\,\,\forall i,\,\,\,0\leq i\leq \left\lfloor \tfrac{d+1}{2}\right\rfloor ,
\end{equation*}
by (\ref{HSCHEN}). On the other hand, (\ref{STAL2}) implies for the
coordinates of $\widetilde{\mathbf{h}}\left( \mathcal{S}\right) $ : 
\begin{equation*}
\widetilde{\frak{h}}_{i}\left( \mathcal{S}\right) =\widetilde{\frak{h}}%
_{\left( d+1\right) -i}\left( \mathcal{S}\right) ,\,\,\,\forall
i,\,\,\,0\leq i\leq d+1,
\end{equation*}
and therefore 
\begin{equation}
\widetilde{\frak{h}}_{i}\left( \mathcal{S}\right) \leq \frak{h}_{i}\left( 
\text{CycP}_{d+1}\left( \frak{b}\right) \right) ,\,\,\,\forall i,\,\,\,0\leq
i\leq d+1\ .  \label{ANISO}
\end{equation}
Hence,{\small 
\begin{eqnarray*}
\frak{f}_{i}\left( \mathcal{S}\right) &=&\sum_{j=0}^{i+1}\ \tbinom{d+1-j}{d-i%
}\ \frak{h}_{j}\left( \mathcal{S}\right) \\
&&\ \medskip \\
&=&\sum_{j=0}^{i+1}\ \tbinom{d+1-j}{d-i}\ \widetilde{\frak{h}}_{j}\left( 
\mathcal{S}\right) \ +\sum\limits_{j=\left\lceil \frac{d}{2}\right\rceil
+1}^{i+1\smallskip }\ \tbinom{d+1-j}{d-i}\ \left( \frak{h}_{d-j}\left(
\partial \mathcal{S}\right) -\frak{h}_{d+1-j}\left( \partial \mathcal{S}%
\right) \right) \medskip \\
&&\  \\
&=&\sum_{j=0}^{i+1}\ \tbinom{d+1-j}{d-i}\ \widetilde{\frak{h}}_{j}\left( 
\mathcal{S}\right) +\sum\limits_{j=d-i}^{\left\lfloor \frac{d}{2}%
\right\rfloor \smallskip }\ \tbinom{j}{d-i}\ \left( \frak{h}_{j-1}\left(
\partial \mathcal{S}\right) -\frak{h}_{j}\left( \partial \mathcal{S}\right)
\right) \medskip \\
&&\  \\
&&\text{[ by interchanging\ }\left( d+1\right) -j\text{ \ and\ }j\text{, and
using } \\
&&\text{the Dehn-Sommerville relations for }\mathbf{h}\left( \partial 
\mathcal{S}\right) \text{]\medskip } \\
&&\  \\
&\leq &\sum_{j=0}^{i+1}\ \tbinom{d+1-j}{d-i}\ \frak{h}_{i}\left( \text{CycP}%
_{d+1}\left( \frak{b}\right) \right) + \\
&&+\sum\limits_{j=d-i}^{\left\lfloor \frac{d}{2}\right\rfloor \smallskip }\ 
\tbinom{j}{d-i}\ \left( \frak{h}_{j-1}\left( \partial \mathcal{S}\right) -%
\frak{h}_{j}\left( \partial \mathcal{S}\right) \right) \medskip \ \ \text{[
by (\ref{ANISO}) ]} \\
&&\  \\
&=&\frak{f}_{i}\left( \text{CycP}_{d+1}\left( \frak{b}\right) \right)
-\sum\limits_{j=d-i}^{\left\lfloor \frac{d}{2}\right\rfloor \smallskip }\ 
\tbinom{j}{d-i}\ \left( \frak{h}_{j}\left( \partial \mathcal{S}\right) -%
\frak{h}_{j-1}\left( \partial \mathcal{S}\right) \right) \\
&&\ 
\end{eqnarray*}
}for all $i,0\leq i\leq d$.\hfill $\square $

\begin{corollary}
\label{SECUBB}Let $\mathcal{S}$ be a simplicial $d$-$\mathbf{ball}$ with $%
\frak{f}_{0}\left( \mathcal{S}\right) =\frak{b}$ vertices. Suppose that $%
\frak{f}_{0}\left( \partial \mathcal{S}\right) =\frak{b}^{\prime }$. Then%
\emph{:\medskip } 
\begin{equation}
\frak{f}_{d}\left( \mathcal{S}\right) \leq \frak{f}_{d}\left( \emph{CycP}%
_{d+1}\left( \frak{b}\right) \right) -\left( \frak{b}^{\prime }-d\right)
\label{UP-LO2}
\end{equation}
\end{corollary}

\noindent \textit{Proof. } For $i=d$, (\ref{WEAKUBTSB}) gives{\small 
\begin{eqnarray*}
\frak{f}_{d}\left( \mathcal{S}\right) &\leq &\frak{f}_{d}\left( \text{CycP}%
_{d+1}\left( \frak{b}\right) \right) -\sum\limits_{j=0}^{\left\lfloor \frac{d%
}{2}\right\rfloor \smallskip }\ \left( \frak{h}_{j}\left( \partial \mathcal{S%
}\right) -\frak{h}_{j-1}\left( \partial \mathcal{S}\right) \right) \\
&&\  \\
&=&\frak{f}_{d}\left( \text{CycP}_{d+1}\left( \frak{b}\right) \right) -\frak{%
h}_{\left\lfloor \frac{d}{2}\right\rfloor }\left( \partial \mathcal{S}\right)
\\
&&\smallskip \  \\
&\leq &\frak{f}_{d}\left( \text{CycP}_{d+1}\left( \frak{b}\right) \right) -%
\frak{h}_{1}\left( \partial \mathcal{S}\right) ,
\end{eqnarray*}
}where the latter inequality comes from the LBT \ref{LBTH} for the
simplicial sphere $\partial \mathcal{S}$. Now obviously, $\frak{h}_{1}\left(
\partial \mathcal{S}\right) =\frak{b}^{\prime }-d$.\hfill $\square $

\section{Crepant Resolutions and UBT\label{SECTION4}}

\noindent {}Let $\left( U_{\sigma },\text{orb}\left( \sigma \right) \right) $
be a Gorenstein toric singularity as in \S \ref{SECTION1} (cf. Thm. \ref
{GORC}).

\begin{theorem}[Necessary Existence Condition]
\label{MAIN}If $U_{\sigma }$ admits a crepant desingularization, then the
normalized volume of the lattice polytope $P_{\sigma }$ \emph{(}defined\emph{%
\ }in \emph{(\ref{PISIGMA})) }has the following upper bound\emph{\footnote{%
By abuse of notation, we write $\mathbb{Z}^{d}$ instead of $\mathbf{H}%
_{\sigma }\cap \mathbb{Z}^{d+1}$ ($\cong \mathbb{Z}^{d}$)}} 
\begin{equation}
\text{\emph{Vol}}_{\mathbf{norm}}\left( P_{\sigma }\right) \leq \frak{f}%
_{d}\left( \emph{CycP}_{d+1}(\#(P_{\sigma }\cap \mathbb{Z}^{d}))\right)
-\left( \#(\partial P_{\sigma }\cap \mathbb{Z}^{d})-d\right)
\label{UB-FINAL}
\end{equation}
\end{theorem}

\noindent {}\textit{Proof}. If $U_{\sigma }$ admits a crepant
desingularization, then there must be a basic triangulation, say $\mathcal{T}
$ of $P_{\sigma }$. Since this $\mathcal{T}$ \ is, in particular, maximal,
we have 
\begin{equation}
\text{vert}\left( \mathcal{T}\right) =P_{\sigma }\cap \mathbb{Z}^{d},\ \ 
\text{ vert}\left( \partial \mathcal{T}\right) =\partial P_{\sigma }\cap 
\mathbb{Z}^{d}.  \label{VERT-EQ}
\end{equation}
On the other hand, 
\begin{equation}
\text{Vol}_{\mathbf{norm}}\left( P_{\sigma }\right) = \frak{f}_{d}\left( \mathcal{T}\right) .  \label{FSX}
\end{equation}
Finally, since $\mathcal{T}$ is a simplicial $d$-ball, one deduces (\ref
{UB-FINAL}) from (\ref{UP-LO2}), (\ref{VERT-EQ}), (\ref{FSX}).\hfill $%
\square $

\begin{example}
Let 
\begin{equation*}
\sigma =\mathbb{R}_{\geq 0}\,e_{1}+\mathbb{R}_{\geq 0}\,e_{2}+\mathbb{R}%
_{\geq 0}\,e_{3}+\mathbb{R}_{\geq 0}\,\left( -3,-7,-9,20\right) \subset 
\mathbb{R}^{4}
\end{equation*}
be the four-dimensional cone supporting the lattice $3$-simplex 
\begin{equation*}
\mathbf{s}_{\sigma }=\text{conv}\left( \left\{ e_{1},e_{2},e_{3},\left(
-3,-7,-9,20\right) \right\} \right) =\left\{ \mathbf{x}\in \sigma \ \left| \
\left\langle \mathbf{m}_{\sigma },\mathbf{x}\right\rangle =1\right. \right\}
,
\end{equation*}
where $\mathbf{m}_{\sigma }=(1,1,1,1).$ Obviously, Vol$_{\mathbf{norm}%
}\left( \mathbf{s}_{\sigma }\right) =20.$ On the other hand, since $\mathbf{s%
}_{\sigma }$ has $8$ lattice points (its $4$ vertices and further $4$
integer points lying in its relative interior), computing the right-hand
side of (\ref{UB-FINAL}), we obtain by (\ref{FDCYCLIC}): 
\begin{equation*}
\tbinom{8-\left\lceil \frac{4}{2}\right\rceil }{\left\lfloor \frac{4}{2}%
\right\rfloor }+\tbinom{8-1-\left\lceil \frac{3}{2}\right\rceil }{%
\left\lfloor \frac{3}{2}\right\rfloor }-(4-3)=\binom{6}{2}+\binom{5}{1}-1=19.
\end{equation*}
This means that $\left( U_{\sigma },\text{orb}\left( \sigma \right) \right) $
necessarily belongs to the class \textbf{(B) }described in \S \ref{SECTION1}.
\end{example}

\begin{remark}
The upper bound (\ref{UB-FINAL}) will be improved considerably in \cite{DHZ2}%
, in the case in which $P$ is a \textit{simplex}, by using a different
combinatorial-topological technique.
\end{remark}

\begin{small}

\end{small}
\end{document}